\documentclass{amsart}
\usepackage[T1]{fontenc}
\usepackage[english]{babel}
\usepackage{amsfonts,amsmath,amssymb,amsthm}  
\usepackage{graphicx,mathtools}
\usepackage{tikz}
\usetikzlibrary{backgrounds,calc}

\input xy
\xyoption{all}

\newcommand{\CC}{\mathbb{C}}

\newcommand{\NN}{\mathbb{N}}

\newcommand{\G}{\operatorname{G}}
\renewcommand{\H}{\operatorname{H}}
\newcommand{\K}{\operatorname{K}}
\renewcommand{\L}{L}
\newcommand{\N}{\operatorname{N}}
\renewcommand{\O}{\operatorname{O}}

\newcommand{\SO}{\operatorname{SO}}

\newcommand{\SU}{\operatorname{SU}}

\newcommand{\U}{\operatorname{U}}
\newcommand{\W}{\operatorname{W}}

\newtheorem{thm}[equation]{Theorem}
\newtheorem{cor}[equation]{Corollary}
\newtheorem{lem}[equation]{Lemma}

\newtheorem*{thm*}{Theorem}
\newtheorem*{cor*}{Corollary}
\newtheorem*{lem*}{Lemma}
\newtheorem*{prop*}{Proposition}
\theoremstyle{definition}
\newtheorem{defn}[equation]{Definition}
\newtheorem*{defn*}{Definition}
\theoremstyle{remark}
\newtheorem{rem}[equation]{Remark}

\newtheorem*{rem*}{Remark}
\newtheorem*{ex*}{Example}
\newtheorem*{assumption*}{Assumption}

\newtheorem*{pf}{Proof}

\newcommand{\cref}[1]{Corollary~\ref{#1}}
\newcommand{\dref}[1]{Definition~\ref{#1}}

\newcommand{\lref}[1]{Lemma~\ref{#1}}

\newcommand{\rref}[1]{Remark~\ref{#1}}

\newcommand{\tref}[1]{Theorem~\ref{#1}}

\numberwithin{equation}{section}

\sffamily

\title{A Note On Polar Representations}

\author[F.~J.~Gozzi]{Francisco J. Gozzi}
\address{CMCC, Universidade Federal do ABC, Av. dos Estados 5001, Santo André - SP, 09210-580, Brazil}
\email{fj.gozzi@gmail.com}
\date{}
\thanks{The author was partially supported by
  FAPESP fellowship 2014/22568-1.}
\subjclass[2010]{58D19, 57S15, 22E46}

\begin{document}
\maketitle

Lie group actions on manifolds are of fundamental importance in geometry and related areas as, for example, in describing the continuous symmetries of a physical system. 
It is then natural that the linear case receives our attention through the study of compact Lie group representations, which are the first examples on the matter and, also, because they provide local models for said, smooth and proper, manifold symmetries. 

We shall discuss a recognition problem in the case of \textit{polar} representations, referring to linear actions of compact Lie groups on euclidean space $V$ admitting a cross-section, i.e., a subspace $\Sigma\subseteq V$ that intersects every orbit orthogonally. 
This is a class of actions that has received special attention from a geometric point of view, see \cite{PalaisTerng} or \cite{BerndtConsoleOlmos}. In particular, they are quotient-equivalent to finite linear group representations, the simplest possible reductions in the copolarity hierarchy of \cite{Gorodski-Lytchak-OnOrbitSpaces}.
Moreover, the geometry of polar representations of connected Lie groups is well understood, they are orbit-equivalent to symmetric space isotropy representations, s-representations for short, after work of J. Dadok \cite{Dadok}. 
In the irreducible case, polar representations have been completely classified up to linear equivalence so that, besides s-representations, they amount to a finite number of exceptions listed in \cite{EH_classif_polarrep}. 
Moreover, an interesting observation by K. Grove and W. Ziller in \cite{GroveZiller} shows that in order to identify an irreducible polar representation it suffices to know the group and a principal isotropy subgroup, $\H\subseteq \G$. 

The aim of the present work is to generalize the previous criterion to encompass general reducible polar representations. 
The key is the existence of a strict fundamental domain, $C\subset \Sigma$, and, more precisely, of the $\G$-isotropy subgroups that occur along its connected strata. This poset of isotropy subgroups 
was already considered in \cite{GroveZiller} as part of \textit{Coxeter polar data} and, in the present work, shall be referred to as a \textit{history} for the representation.
This leads to our main result.
\begin{thm*}
A polar representation of a connected compact Lie group $\G$ is determined up to linear isomorphism by a history and dimension. 
\end{thm*}
The same statement holds true for a disconnected group provided we enhance the hypothesis on the action to be a \textit{Coxeter polar} representation, see \cref{coro Coxeter}.

\smallskip
In the following, after duly introducing some background, we provide a direct proof of our theorem. In the second section, we treat the case of Coxeter polar actions, extending our results to disconnected groups and discussing its application to equivariant reconstruction of polar $\G$-manifolds.

\section{Polar Representations} 
The following brief exposition on the geometry of polar representations of compact connected Lie groups is mostly self-contained,  
with the crucial exception of some of the aforementioned results. We specially rely on the orbit-equivalence to s-representations and their associated Weyl groups, references are \cite{Dadok,PalaisTerng}. We refer to \cite{GroveZiller} for a detailed account of most of the following constructions. For background on finite reflection groups the reader may refer to \cite{bourbaki-reflection}. 

\medskip 
Let $\G$ be a compact connected Lie group and $\rho:\G\to \O(V)$ be an orthogonal representation admitting a linear subspace $\Sigma\subset V$
that intersects every orbit and does so orthogonally, i.e., a polar representation. 
The \textit{polar section} $\Sigma$ integrates the orthogonal distribution to the orbit foliation, when restricted to the open dense set of points whose orbits have maximum dimension. Hence, generically, we have:
\begin{equation}\label{V split}
T_p(\G \cdot p) \oplus \Sigma = V.     
\end{equation}
The existence of a section helps us understand the geometry of the orbit space by optimally reducing the dimensionality of the problem. To be precise, consider the normalizer subgroup, $\N(\Sigma)=\{g\in \G: ~g\cdot \Sigma = \Sigma \}$ and its natural action on the section $\Sigma$. This is a discrete action whose kernel is given by the point-wise stabilizer or centralizer subgroup of $\Sigma$ 
which, in turn, agrees with the -only- principal isotropy subgroup $\H\subset \G$ occurring along the section. 
The quotient $\W:=\N(\Sigma)/\H$ is the \textit{polar group} of the representation, a finite group whose effective action on $\Sigma$ gives the same quotient as the original representation, so that:
\begin{equation} \label{polar group on section}
\Sigma/\W \cong V/ \G.     
\end{equation}
Already from a metric point of view, this quotient space displays an orbifold structure and, hence, encodes the linear representation of the finite group $\W$ on $\Sigma$. Moreover, given that polar representations are orbit-equivalent to s-representations, it follows from the latter that $\W$ acts as a reflection group on the section.
As a classical example, recall that the adjoint representation of a compact connected Lie group on its Lie algebra is polar with section given by a Cartan subalgebra and $\W$ its Weyl group. 

In general, once the orbit space is described as a quotient by a finite group it easily follows that we may find a fundamental domain. 
In the case of a reflection group the situation is even better, as we may find a strict fundamental domain $C$, known as a \textit{chamber}.
In order to do so, consider the complement in $\Sigma$ to the hyperplanes supporting the reflections in $\W$, then choose any connected component and take its closure, this is $C$. We have:  
\begin{equation}\label{chamber and quotient}
C \cong \Sigma/\W \cong V/ \G.
\end{equation}
It is then clear that the orbit space is a Coxeter orbifold, i.e., one whose local orbifold groups are finite linear reflection groups (see \cite{reflection-manifolds}, c.f. Coxeter domain in \cite{GroveZiller}). 

The discreteness of the polar group -an algebraic counterpart to the existence of a section- together with \eqref{chamber and quotient} imply that the isotropy subgroups occurring along connected strata of $C\subseteq \Sigma$ are constant. This motivates the next concept. 
\begin{defn}\label{defn history}
An admissible \textit{history} for a polar representation of a compact connected Lie group $\G$ is the collection $\{\K_I\}_I$ of $\G$-isotropy subgroups occurring along a given chamber $C$.  
\end{defn} 
Observe that the isotropy groups in the history together with the given inclusion relations form a finite lattice, 
where the group $\G$ is the maximum and the principal isotropy group $\H$ along the section is the minimum.

As a general remark, in the context of proper Lie group actions on manifolds, recall that the points in the orbit space
can be stratified by different criteria: their orbit isotropy type, the metric structure of their tangent unit sphere inside $V/\G$ or the isomorphism class of the slice representation to the orbit.
When considering connected strata these three criteria agree. 

The main advantage in dealing with this kind of polar representations is that there is a strict fundamental domain $C\subset \Sigma$ for both the $\G$-action on $V$ and that of its polar group $\W$ on the section $\Sigma$. In particular, strata of $C$ are associated to respective $\G$ and $\W$-isotropy subgroups, which are in $1-1$ correspondence with each other.

Notice that the $\G$-isotropy subgroups $\K_j$ for the open codimension one faces of $C$ are the next to minimal groups $\H\subset \K_j \subset \G $ in the given history for $\rho$ and, in turn, correspond to reflections $r_j\in \W$ that generate the polar group. 
Moreover, each element $r_j$ can be identified as the unique involution in  $\N_{\K_j}(\H)/\H$, the normalizer of the principal isotropy group $H$ inside $\K_j$. For this, observe that the group  $\N_{\K_j}(\H)/\H$ acts transitively on a the sphere normal to the stratum and with trivial principal isotropy, so that $\N_{\K_j}(\H)/\H=\mathbb{S}^i$, $i=0,1,3$.
As a consequence, knowing just the next to minimal $\G$-isotropy subgroups, i.e., the first "layer" of a history for $\rho:\G\to \O(V)$, we can find these generating reflections and, thus, determine $\W$ as a subgroup inside $\N(\H)/\H$. These reflections $S= \{r_j\}_j$ determine a Coxeter system for the action of $\W$ on $\Sigma$ since they are precisely those reflections supported on the hyperplane faces of a chamber. 
We have thus proved:
\begin{lem}\label{lem1 Coxeter system}
A history $\{\K_I\}_I$ for a polar representation $\rho:\G\to \SO(V)$ of a compact connected Lie group $\G$ determines the associated polar group $\W$ and a Coxeter system $(\W, S)$ for its action as a reflection group on the section $\Sigma$.
\end{lem}

We may now exploit the rigid structure of finite Coxeter groups, namely, their classification via Dynkin diagrams. For example, the number of components the asociated Dynkin diagram corresponds to the number of irreducible factors $l\in \NN$ into which $(\W,S)$ splits as a reflection group.
Hence, 
\begin{equation}\label{eq split polar group}
W=\W_1 \times \cdots \times \W_l,
\end{equation} 
where each factor is a reflection group on its own, acting irreducibly on a euclidean space $\Sigma_i$, whose dimension may also be extracted from the corresponding Coxeter system.
In this way, the linear section $\Sigma$ -where $\W$ acts- splits as a product of the previous factors $\Sigma_i$, $i>0$: 
\begin{equation}\label{eq split section}
Sigma= \Sigma_0 \oplus \Sigma_1 \oplus \cdots \oplus \Sigma_l,
\end{equation}
plus, possibly, a trivial summand $\Sigma_0$ to agree with the original polar section. In particular, the dimension of $\Sigma_0$ is just the difference of the cohomogeneity of the action and the sum of the other, known, factors  $\Sigma_i$, $i>0$. Recall that, 
$$ cohom(\G,V)= n - dim(\G) + dim(\H), $$
where $n$ is the dimension of the representation and $\H$ is a principal isotropy group, the minimum in the given history. We have a corollary to the previous lemma. 
\begin{cor}\label{cor history and C}
A history $\{\K_I\}_I$ for a polar representation $\rho:\G\to \SO(V)$ of a compact connected Lie group $\G$ determines the structure of the quotient $V/\G$. 
To be precise, 
the quotient is isometric to a product:
$$ V/ \G = \Sigma_0 \times (\Sigma_1/\W_1) \times \cdots \times (\Sigma_l/\W_l),$$
where each $\W_i$ acts as an irreducible reflection group on $\Sigma_i$, $i>0$, and $\Sigma_0$ is an euclidean factor.
\end{cor}
For an irreducible representation $\rho:\G\to V$ we have that its associated reflection group and Coxeter system $(\W,S)$ is also irreducible -though not necessarily as a group per se. In the irreducible case, there is available a stronger result than that of our Theorem. 
\begin{lem}\label{lemma irreducible recognition GZ} (from \cite{GroveZiller})
An irreducible polar representation of a compact connected Lie group is linearly determined by its dimension, the group, and a principal isotropy subgroup. 
\end{lem}
\begin{pf}
This Lemma is observed in passing in the original work and follows from a direct case by case analysis of the representations, whose classification is well-known. Essentially, a single Lie group can have but a handful of irreducible representations that are polar and it turns out that the principal isotropy type is enough to distinguish among them.
Moreover, the case of s-representations can be found in \cite{tamaru}, computed on the Lie algebra level. The few exceptions, i.e., the irreducible polar representations of compact connected Lie groups which are not themselves s-representations proper, are in fact orbit-equivalent subgroups of s-representations (see \cite{EH_classif_polarrep}). 
In this way, if one knows the explicit inclusion of the original group $\L$ into a larger one that gives the corresponding s-representation, $\Tilde{\L}$, then the principal isotropy for $\L$ can be determined by simply restricting that of $\Tilde{\L}$.  \qed
\end{pf}
Towards the end of the present work we provide an alternative proof of our Theorem which does not rely on this last lemma, see \rref{rem alternative proof}. 

\medskip
\noindent\textbf{Proof of the Theorem:}\\
We consider a polar representation $\rho: \G\to V$ of a compact connected Lie group as our incognitum, of which we are only given a lattice of subgroups, $\{\K_I\}_I$, corresponding to an admissible history for the representation. 

Firstly, as shown in \lref{lem1 Coxeter system}, we determine the associated polar group and its action on the section. The factorization of the Coxeter system $(\W,S)$ into irreducible components forces the original representation $\rho$ to split like-wise. In particular, keeping with previous notations from \lref{lem1 Coxeter system} and \cref{cor history and C}, we have that $\rho:\G\to \SO(V)$ factors as: 
$$ V=V_0 \oplus V_1\oplus \cdots \oplus V_l.$$
This is a direct sum of $l\in\NN$ non-trivial irreducible summands  $V_i$ plus, possibly, a subspace $V_0=\Sigma_0$ where the $\G$-action is trivial. 
This observation should be clear
once the Coxeter system $(\W,S)$ is an invariant for the class of orbit-equivalent representations, in fact, for quotient-equivalent ones. 

Back to the splitting of the section $\Sigma$ for the action of $\W$, notice that each factor $\Sigma_i$ in \eqref{eq split section} is point-wise fixed by the action of the group:
$$\W_{\Sigma_i}:= (\W_1 \times \cdots \times \hat{\W_i} \times \cdots \W_l) \subseteq \N(\H)/\H.$$   
Then, the $\G$-isotropy group of a generic (regular) point $p_i \in \Sigma_i $ satisfies
\begin{equation}
\G_{p_i} \supseteq \W_{\Sigma_i}\cdot \H. 
\end{equation}
In fact, $\G_{p_i}$ is the minimal group in the history lattice satisfying the above inclusion, thanks to the correspondence between the history and the $\W$-isotropies (once $p_i$ can be chosen to lie in $C$).
Putting all together for the i-th factor, we have that the group $\G$ acts on $V_i$ admitting principal isotropy group $\G_{p_i}= \K_{I_i}$. 

Finally, the latter factor representation can be made effective by standard means. For this take the quotient of $\G$ by the maximal group inside a principal isotropy group, $\N_i \subseteq \K_{I_i}$, that is normal under $\G$. 
The outcome effective and irreducible subrepresentation, 
$\rho_i:\G/\N_i \to \SO(V_i)$, with principal isotropy $ \K_{I_i}/\N_i$,
can be identified using \lref{lemma irreducible recognition GZ}. 

The original representation $\rho:\G\to \SO(V)$ can be retrieved as follows: 
$$\G \xrightarrow{~\Delta~} \G/\N_1 \times \cdots \times \G/\N_l \xrightarrow{~\oplus_i \rho_i~} \SO(V_0\oplus V_1\oplus \cdots \oplus V_l).$$
The dimension of $V$ is only required to determine the fixed subspace $V_0$.
This concludes the proof of the Theorem. 
\qed

Notice that if we were given all the $\G$-isotropy subgroups along a section instead of precisely a chamber, the same proof would work, after suitably extracting a Coxeter system out of more reflections for $\W$.

\section{Polar $\G$-manifolds}

This section surveys the main application of our Theorem into equivariant reconstruction of polar $\G$-manifolds. 
\smallskip

A proper isometric action of a Lie group $\G$ on a complete riemannian manifold $M$ is said to be polar provided there is an analogous cross-section $\Sigma$ which, in this context, corresponds to a complete immersed submanifold intersecting every $\G$-orbit orthogonally. 
The polar section $\Sigma$ and its $\G$-translates are integral leaves to the totally geodesic distribution given as the orthogonal complements to the orbits, if we restrict to the open dense set of regular points.

The analogous definition of the polar group as $\W= \N(\Sigma)/\H$ grants a discrete group acting in a proper and isometric manner on $\Sigma$.
In particular, \eqref{polar group on section} holds, if $V$ is suitably replaced for the ambient manifold $M$.
In this context, a reflecting element is understood as an involution $r\in W$ that fixes a hypersurface. It is then possible to construct a domain $C$, a chamber, by an analogous procedure to that of linear reflection groups, namely, by taking the closure to a connected component of the complement in $\Sigma$ to all supporting hypersurfaces fixed by reflections in $\W$. 
In general, this chamber $C$ does not provide a fundamental domain for $M$, but motivates the following.
\begin{defn} 
A polar action is \textit{Coxeter polar} if its polar group $\W$ is generated by singular reflections and the chamber $C$ has trivial normalizer inside $\W$. 
\end{defn}
Here, singular reflection stands for a reflecting element of the $\W$-isotropy that comes from a singular $\G$-isotropy subgroup. 

 Notice that if two points in $C$ lie over the same orbit, then there is an element in $\W$ that takes one to the other and preserves $C$, since the $\W$ isotropy at a point $p\in C\subseteq \Sigma$ acts transitively on the set of chambers that contain the point. Hence, the latter condition in the above definition implies that the chamber $C$ is a strict fundamental domain. Again, mutatis mutandi, \eqref{chamber and quotient} holds. 

Back to the linear case, observe that polar representations by compact connected Lie groups are always Coxeter polar, as follows from the following lemma.
\begin{lem}\label{lemma 2.4 GZ}  A polar representation $(\G,V)$ is Coxeter polar if and only if the action of the group is orbit equivalent to that of its identity component $(\G_0,V)$. 
\end{lem}
\begin{pf}
This is item c) of Lemma 2.4 in \cite{GroveZiller}.\qed 
\end{pf}
Notice that the notion of a history extends naturally to the case of Coxeter polar representation (c.f. \dref{defn history}).
We have the following generalization of our Theorem for disconnected $\G$.
\begin{cor}\label{coro Coxeter}  
A Coxeter polar representation of a compact Lie group $\G$ is determined by a history and dimension.
\end{cor}

\begin{pf}
Given a history $\{K_I\}_I$ for a representation $(\G,V)$ we can restrict the subgroups to the identity component and, thus, obtain a history $\{K_I\cap \G_0\}_I$ for $\G_0$. It follows from our main Theorem that we can determine the Coxeter polar representation of $\G_0$ on $V$.
Then, making use of the previous Lemma, we have that for $\G$ to be orbit equivalent to $\G_0$, principal orbits need to coincide, so that: 
$$\G/\H = \G_0 / (\G_0 \cap \H).$$
This implies that $\G$ is generated by $\G_0$ and the principal isotropy group $\H$. We need only determine the action of $\H$. Observe that $\H$ acts trivially over the linear section $\Sigma$ and by the isotropy representation along the tangent space to a principal orbit, i.e., the $\H$-action on $T_{\H}(\G\!/\!\H)$. Here we are using \eqref{V split}. 
\qed
\end{pf}
Our \cref{coro Coxeter} is a linear recognition tool designed to fit
into the general framework of polar equivariant reconstruction of $\G$-manifolds, specifically, as a technical step in the proof of \tref{thm GroveZiller}. We refer to the original work in \cite{GroveZiller} for an in depth treatment. 

The same proof as above, allows for the generalization of \lref{lemma irreducible recognition GZ} to the case of a disconnected group with a Coxeter polar representation. 
\begin{cor}\label{coro Coxeter irreducible} 
An irreducible Coxeter polar representation of a compact Lie group $\G$ is determined by its dimension, the group and a principal isotropy subgroup.
\end{cor}
Without the the irreducibility assumption this would be the second part of item b) in Lemma 2.4 from \cite{GroveZiller}, which is flawed as  the following counter-example shows. We propose that the technical gap left by this issue be replaced by our \cref{coro Coxeter}. 
\begin{rem}\label{counter-example}
Two polar representations of the same group and dimension may even have coincident principal isotropy groups without being isomorphic. 
For this, we need only compare a factor-wise product action of $\SU(n)\times\U(1)$ on $\CC^n\times\CC$ with the representation given as a cover of the standard $\U(n)$ action on $\CC^n$ times the circle action on $\CC$. Both have the same standard $\SU(n-1)$ subgroup as a principal isotropy subgroup. 
\end{rem}

\begin{rem}\label{rem Coxeter polar marking}
A (Coxeter) polar $\G$-manifold $M$ has (Coxeter) polar slice representations. 
In fact, the tangent space to a polar section $\Sigma \subset M$ passing through a point $p$ is a linear section for the corresponding representation of the isotropy group $\G_p$ on the slice $T_p (\G\cdot p)^{\perp}$ (we refer to \cite{HLO_isoparametric_submanifolds} for the polar case and to \cite{GroveZiller} for the enhancement to Coxeter polar). 
It follows that, the isotropy subgroups along the strata of a chamber $C\subseteq \Sigma$ correspond, locally near a point $p\in C$, to a history for the corresponding Coxeter polar slice representation of $\G_p$ on $T_p(\G\cdot p)$. 
\end{rem}
This leads to the following notion. 
\begin{defn}\label{defn coxeter polar data}
 \textit{Coxeter polar data} for a $\G$-manifold $M$ is defined as a Coxeter orbifold $C$ 
 together with a \textit{compatible} marking of its strata by (constant) $\G$-subgroups, $\G(C):=\{\G_c\}_{[c]\in Str(C)}$. 
Compatibility for the pair $(C,\G(C))$  means that the marking on $C$ is locally (near $p \in C$) that of a history for a Coxeter polar representation, whose polar group agrees with the corresponding local orbifold group at $p$. 
 \end{defn}

\begin{rem}\label{rem equivariant homeo}   
In general, knowing a fundamental domain for an action together with the point-wise isotropy subgroups allows to reconstruct each orbit separately. 
Hence, from a simplified point of view, a Coxeter polar $\G$-manifold with data $(C, \G(C))$ can be determined up to equivariant homeomorphism as a quotient of $\G\times C$ in the following way (cf. \cite{HambletonHausmann}):
\begin{equation}\label{eqn equivariant homeo for fundamental domain}
M= \G\times C/\!\sim ~~\text{with:}~~ (g,c)\!\sim\! (g',c')\iff\mbox{c=c'}~\wedge~ g^{-1}g'\in \G_c.
\end{equation}
The obvious problem is that we loose track of the smooth structure. 
\end{rem}

\begin{rem}\label{rem alternative proof}
Back to the linear case, this suggests an alternative proof of our main Theorem. The argument should begin in the same fashion in order to determine the structure of the quotient. Recall \cref{cor history and C}, which says that the algebraic part of Coxeter polar data for a linear representation, i.e., its history, determines the whole Coxeter polar data. 
Then, one can apply \eqref{eqn equivariant homeo for fundamental domain}
to retrieve the $\G$-equivariant homeomorphism type of the action.
Finally, a strong rigidity from result \cite{LeeWasserman} states that for a compact connected Lie group representation topological equivalence implies linear equivalence. 
\end{rem}

\begin{rem} In the case of a proper isometric action of a Lie group on a manifold, the semi-local structure of the foliation and the smoothness issues pointed out in \rref{rem equivariant homeo} can be addressed via the classical Slice Theorem which is, in essence, an equivariant reconstruction tool. 
For this, 
assume that the slice representation at a point $p\in M$ is known, say $\G_p\to \SO(V)$, with $V= T_p (\G\cdot p)^{\perp}$. Then, an $\epsilon$-neighbourhood, $B_{\epsilon}(\G\cdot p)$, around this orbit can be described up to equivariant diffeomorphism as follows:
$$ B_{\epsilon}(\G\cdot p)~~ \cong~~ {\G \times_{\G_p} V}. $$
The term on the right is the quotient of $\G \times V$ by the diagonal action of ${\G_p}$, as a subgroup of $\G$ from the right and by the previous representation on $V$ from the left (see \cite{Bredon}). 
\end{rem}
It then follows that we can perform a semi-local reconstruction of a Coxeter polar $\G$-manifold $M$ for which there is available Coxeter polar data $(C,\G(C))$. 
Our \cref{coro Coxeter} serves as a technical step in determining the corresponding slice representations. It is in this sense that our work serves the proof of the following stronger result.
\begin{thm}\label{thm GroveZiller}(from \cite{GroveZiller}) 
There is a canonical construction of a Coxeter polar $\G$-manifold $M(D)$ for Coxeter polar data $D = (C, \G(C))$. 
Furthermore, if $D$ is Coxeter polar data for a Coxeter polar $\G$-manifold $M$, then $M(D)$ is equivariantly diffeomorphic to $M$. 
\end{thm}

\medskip
The issue of the current paper arouse during the author's Ph.D. studies at IMPA, thanks go to his former advisors Luis A. Florit and W. Ziller.


\end{document}